\documentclass[12pt]{amsart}

\newcommand{\Q}{{\mathbf Q}}

\newcommand{\Kbar}{{\overline{K}}}

\newcommand{\Dbar}{{\overline{D}}}

\newcommand{\Fbar}{{\overline{\F}}}
\newcommand{\Z}{{\mathbf Z}}

\newcommand{\TT}{{\mathcal T}}
\newcommand{\OO}{{\mathcal O}}
\newcommand{\F}{{\mathbf F}}
\newcommand{\mm}{{\mathfrak m}}

\newcommand{\Gal}{\operatorname{Gal}}

\newcommand{\nichts}{{\left.\right.}}
\newcommand{\isom}{\cong}
\newcommand{\tensor}{\otimes}

\newtheorem{theorem}{Theorem}
\newtheorem{lemma}[theorem]{Lemma}

\newtheorem{prop}[theorem]{Proposition}

\theoremstyle{definition}

\theoremstyle{remark}
\newtheorem{rem}{Remark$\!\!$}

\begin{document}

\title[Zeros of sparse polynomials]{Zeros of sparse polynomials
	over local fields of characteristic~$p$}
\author{Bjorn Poonen}
\address{Department of Mathematics, University of California, Berkeley, CA 94720-3840, USA}
\email{poonen@math.berkeley.edu}
\date{March 2, 1998}


\maketitle

\section{Statement of results}
\label{introduction}

Let $K$ be a field of characteristic~$p>0$
equipped with a valuation $v: K^\ast \rightarrow G$
taking values in an ordered abelian group $G$.
Let $\OO_K=\{\alpha \in K: v(\alpha) \ge 0\}$
and $\mm_K=\{\alpha \in K: v(\alpha) > 0\}$
be the valuation ring and maximal ideal, respectively,
and suppose that the residue field $\OO_K/\mm_K$
is finite, with $q$ elements.

\begin{theorem}
\label{maintheorem}
If $f(x) = a_0 x^{n_0} + a_1 x^{n_1} + \cdots + a_k x^{n_k}$
is a polynomial with $k+1$ nonzero coefficients $a_i \in K^\ast$,
then $f$ has at most $q^k$ distinct zeros in $K$.
\end{theorem}

This upper bound is sharp:
if $K$ is $\F_q((T))$ with the usual
discrete valuation $v:K^\ast \rightarrow \Z$,
if $V \subset K$ is an $\F_q$-subspace of dimension $k$,
and if $c \in K$ is nonzero,
then the polynomial $f(x) := c \prod_{\alpha \in V} (x-\alpha)$
has the form $a_0 x + a_1 x^q + \cdots + a_k x^{q^k}$
for some $a_0, a_1, \ldots, a_k \in K^\ast$.

Theorem~\ref{maintheorem} is the case $d=1$
of the following generalization, which bounds the number
of distinct zeros of {\em bounded degree}.
Let $\mu(n)$ be the M\"obius $\mu$-function.

\begin{theorem}
\label{degreetheorem}
Fix $d \ge 1$.
If $f(x) = a_0 x^{n_0} + a_1 x^{n_1} + \cdots + a_k x^{n_k}$
is a polynomial with $k+1$ nonzero coefficients $a_i \in K^\ast$,
then the number of distinct zeros of $f$ in $\Kbar$ of degree
at most $d$ over $K$ is at most $\sum_{j=1}^d \sum_{i|j} q^{ik} \mu(j/i)$.
\end{theorem}

This upper bound is sharp as well, for every $q$, $k$, and $d$.
Let $K=\F_q((T))$ and $v$ be as before.
Let $\F$ be a finite field containing $\F_{q^i}$ for $i \le d$.
Let $V \subset \F((T))$ be a $k$-dimensional $\F$-vector space
that is $\Gal(\F/\F_q)$-stable (or equivalently, has an $\F$-basis
of elements of $K$).
Then equality is attained in Theorem~\ref{degreetheorem}
for $f(x) := c \prod_{\alpha \in V} (x-\alpha)$
for any $c \in K^\ast$.
(The inner sum in Theorem~\ref{degreetheorem} performs
the inclusion-exclusion to count zeros of {\em exact} degree $j$.)

We make no claim that these are the only polynomials that attain
equality; in fact there are many others.
For example, if $K$, $V$, and $f$ are as in the previous paragraph, and if
the $\F$-basis of $V$ consists of elements of $K$ of distinct valuation,
with all these valuations divisible by a single integer $e \ge 1$,
then $f(x^e)$ also attains equality,
as a short argument involving Hensel's lemma shows.
Other examples can be constructed using the observation
that if $f(x) \in K[x]$ has $N$ zeros in a given field extension $L$ of $K$,
one of which is 0,
then the same holds for $x^m f(1/x)$ when $m > \deg f$.

\begin{rem}
H.~W.~Lenstra, Jr.~\cite{lenstra}
proves related facts for finite extensions $L$ of $\Q_p$,
using very different methods.
One of his results is that for any such $L$ and any positive integer $k$,
there exists a positive integer $B=B(k,L)$ with the following property:
if $f \in L[x]$ is a nonzero polynomial with at most $k+1$ nonzero
terms and $f(0) \not= 0$, then $f$ has at most $B$ zeros in $L$,
{\em counted with multiplicities}.
His bound $B(k,L)$ is explicit, but almost certainly not sharp.
Finding a sharp bound seems difficult in general, although Lenstra
does this for the case $k=2$ and $L=\Q_2$ (the bound then is 6).
He also applies his local result to bound uniformly
the number of factors of given degree over number fields.
In~\cite{lenstra2} he shows that if $f$ is represented sparsely,
then these factors can be found in polynomial time.
\end{rem}

\begin{rem}
We cannot count multiplicities in either of {\em our} theorems and hope to
obtain a bound depending only on $k$ and $K$
(and $d$, for Theorem~\ref{degreetheorem}),
because of examples like $f(x)=(1+x)^{q^m}$ with $m \rightarrow \infty$.
Requiring that $f$ not be a $p$-th power would not eliminate the problem,
because one could also take $f(x)=(1+x)^{q^m+1}$.
\end{rem}

\section{Proof of Theorem~\ref{maintheorem}}
\label{proof1}

By a {\em disk} in a valued field $K$,
we mean either
an ``open disk'' $D(x_0,g):=\{x \in K: v(x-x_0)>g\}$,
or a ``closed disk'' $\Dbar(x_0,g):=\{x \in K: v(x-x_0) \ge g\}$
where $x_0 \in K$ and $g \in G$.

Let $\sigma_1$, $\sigma_2$, \dots, $\sigma_t$ be the
non-vertical segments of the Newton polygon of $f$.
Let $-g_j \in G \tensor \Q$ be the slope of $\sigma_j$.
If $e_1,e_2,\ldots,e_r$ are the exponents
of the monomials in $f$ corresponding to points on a given $\sigma_j$,
define $N_j$ as the largest integer for which
the images of $(1+x)^{e_1}$, $(1+x)^{e_2}$, \dots, $(1+x)^{e_r}$
in $\F_p[x]/(x^{N_j})$ are linearly {\em dependent} over $\F_p$.
We say that the $\sigma_j$ are in a {\em proper order}
if $N_1 \ge N_2 \ge \cdots \ge N_t$.
This particular ordering is crucial to the proof,
but it is hard to motivate its definition.
It was discovered by analyzing proofs
of many special cases of Theorem~\ref{maintheorem}.
For instance, if the Newton polygon of $f$ has $k$ non-vertical segments
(each associated with exactly two exponents),
then the segments are being ordered according to the $p$-adic absolute values
of their horizontal lengths.

\begin{lemma}
\label{distances}
Let $L$ be a field of characteristic~$p>0$
with a valuation $v: L^\ast \rightarrow G$.
Suppose $f(x) = a_0 x^{n_0} + a_1 x^{n_1} + \cdots + a_k x^{n_k} \in L[x]$
with each $a_i$ nonzero.
List the segments of the Newton polygon of $f$ in a proper order as above.
Fix $u$ and let $-g_u \in G \tensor \Q$ be the slope of
the $u$-th segment $\sigma_u$.
Suppose $r \in L$ is {\em not} a zero of $f$, and $v(r)=g_u$.
Let $S$ be the set of zeros of $f$ in $L$ lying inside $D(r,g_u)$.
Then $\# \{v(\alpha-r): \alpha \in S \} \le k+1-u$.
\end{lemma}

\begin{proof}
Replacing $f(x)$ by $cf(rx)$ for suitable $c \in L^\ast$,
we may reduce to the case in which $r=1$, $g_u=0$, $f \in \OO_L[x]$,
and $f \bmod \mm_L$ is nonzero.
Write $f(1+x)=\sum_{j=0}^{n_k} b_j x^j$,
and let $M$ be the smallest integer for which $b_M$ is nonzero modulo $\mm_L$.
By definition of $N_u$, $f(1+x) \not\equiv 0 \bmod (\mm_L,x^{N_u+1})$.
Hence $M \le N_u$.

For each $i \le u$, we have $N_u \le N_i$, so $M \le N_i$,
and there is some $\F_p$-linear relation in $\F_p[x]/(x^M)$
between the $(1+x)^e$ for the exponents $e$ associated to $\sigma_i$.
The subspace of $\F_p[x]/(x^M)$ spanned by the $(1+x)^e$,
where $e$ ranges over {\em all} the exponents in $f$,
then has dimension at most $(k+1)-u$,
since the $u$ relations above are independent,
the largest $e$ involved in each relation being distinct from
the others.
It follows that the $\F_p$-subspace of $K$
spanned by $b_0$, $b_1$, \dots, $b_{M-1}$
is at most $(k+1-u)$-dimensional.
Then $\#\{v(b_i): 0 \le i < M \text{ and } b_i \not=0 \} \le k+1-u$,
because nonzero elements of distinct valuations
are automatically $\F_p$-independent.
The left endpoints of the
negative slope segments of the Newton polygon of $f(1+x)$
correspond to $b_i$ of distinct valuations for $i<M$,
so there are at most $k+1-u$ such segments.
Hence at most $k+1-u$ positive elements of $G$ can
be valuations of zeros of $f(1+x)$,
which is what we needed to prove.
\end{proof}

\begin{rem}
Note that there is no assumption on the residue field
in Lemma~\ref{distances}; $L$ could even be algebraically closed.
\end{rem}

Let $S$ be any finite subset of a field $L$ with valuation $v$.
We associate a tree $\TT$ to $S$ as follows.
(See~\cite{stanley} for arboreal terminology.)
Let $\TT$ be the Hasse diagram of the finite poset
(ordered by inclusion) of nonempty sets
of the form $S \cap D$ where $D$ is a disk.
Clearly $\TT$ is a tree, whose leaves are the singleton subsets of $S$.
We would obtain the same tree if we required the disks $D$ to be
open (resp.\ closed), since $S$ is finite.

Suppose $r$ and $S$ are as in Lemma~\ref{distances}.
Let $T_0 > T_1 > \cdots > T_\ell$ be the longest chain in $\TT$.
Then $T_\ell$ is a leaf, and $\# T_\ell=1$.
Choose $r_0 \in D(r,g_u) \setminus S$ closer to the element of $T_\ell$
than to any other element of $S$.
For various $g>g_u$, the set $S \cap D(r_0,g)$ can equal
$T_0$, $T_1$, \dots, $T_\ell$, or $\emptyset$.
Hence
	$$\# \{v(\alpha-r_0): \alpha \in S\} \ge \ell+1.$$
On the other hand, Lemma~\ref{distances} applied to $r_0$ yields
	$$\# \{v(\alpha-r_0): \alpha \in S\} \le k+1-u.$$
Combining these, we have that {\em the length $\ell=\ell(\TT)$
of the tree satsifies $\ell \le k-u$.}

Suppose $S_0 \in \TT$ is not a leaf (i.e. $\#S_0>1$),
and let $g=\min\{v(s-t):s,t \in S_0\}$,
so that for any $s \in S_0$,
$\Dbar(s,g)$ is the smallest disk containing $S_0$.
Then the children of $S_0$ in the tree
are nonempty sets of the form $S \cap D(x_0,g)$
for some $x_0 \in \Dbar(s,g)$.
In particular the number of children is at most the
size of the residue field of $L$.

\begin{proof}[Proof of Theorem~1]
Let notation be as in Lemma~\ref{distances}, but take $L=K$.
By the theory of Newton polygons, each nonzero zero of $f$
has valuation equal to $g_u$ for some $u$.
Let us now fix $u$ and let $Z_u$
be the number of zeros in $K$ of valuation $g_u$.
We may assume $g_u \in G$, since otherwise $Z_u=0$.
Then $\{x \in K: v(x)=g_u\}$ is the union of $q-1$ open disks $D_j$
of the form $D(x_j,g_u)$.
As above, the tree corresponding to the set of zeros in $D_j$
has length at most $k-u$, and each vertex has at most $q$ children.
Hence the tree has at most $q^{k-u}$ leaves,
and $Z_u \le (q-1)q^{k-u}$.
Allowing for the possibility that $0$ also is a zero of $f$,
we find that the number of zeros of $f$ in $K$ is at most
	$$1+\sum_{u=1}^t Z_u \le 1 + \sum_{u=1}^t (q-1)q^{k-u} \le 1 + \sum_{u=1}^k (q-1)q^{k-u} = q^k.$$
\end{proof}

\section{Valuation theory}

Before proving Theorem~\ref{degreetheorem},
we will need to recall some facts from valuation theory.
We write $(K,v)$ for a field $K$ with a valuation $v$.
We say that $(L,w)$ is an {\em extension} of $(K,v)$
if $K \subseteq L$ and $w|_K=v$.
In this case,
when we say that $L$ has the same value group (resp.\ residue field)
as $K$,
we mean that the inclusion of value groups (resp.\ residue fields)
induced from the inclusion of $(K,v)$ in $(L,w)$
is an isomorphism.
Recall that any valuation on a field $K$ admits at least
one extension to any field containing $K$.
An abelian group $G$ is {\em divisible} if for all $g \in G$ and $n \ge 1$,
the equation $nx=g$ has a solution $x$ in $G$.

\begin{prop}
\label{makedivisible}
Any valued field can be embedded in another valued field
having the same residue field, but divisible value group.
\end{prop}

\begin{proof}
Let $v: K^\ast \rightarrow G$ be the original valuation.
If $G$ is not already divisible, then there exists $g \in G$
and a prime number $n$ such that $n x=g$ has no solution in $G$.
Pick $\alpha \in K^\ast$ with $v(\alpha)=g$, and extend
$v$ to a valuation on $L=K(\alpha^{1/n})$.
Let $e$ and $f$ denote the ramification index and residue class degree
for $L/K$.
Then $e=n$, and the inequality $ef \le n$
(Lemma~18 in Chapter~1 of~\cite{schilling})
forces $f=1$.
An easy Zorn's lemma argument now shows that $v$ extends
to a valuation $v:M^\ast \rightarrow G \tensor \Q$ where
$M$ is an extension with the same residue field as $K$,
but with divisible value group.
\end{proof}
Recall that $(L,w)$ is called an {\em immediate extension} of $(K,v)$ if
\begin{enumerate}
	\item $(L,w)$ is an extension of $(K,v)$;
	\item $(L,w)$ has the same value group as $(K,v)$; and
	\item $(L,w)$ has the same residue field as $(K,v)$.
\end{enumerate}
Also recall that $(K,v)$ is called {\em maximally complete}
if it has no nontrivial immediate extensions.

\begin{prop}
\label{krullresult}
Every valued field has a maximally complete immediate extension.
\end{prop}

\begin{proof}
This is an old result of Krull:
see Theorem~5 of Chapter~2 in~\cite{schilling}.
\end{proof}

\begin{prop}
\label{embedFq}
Suppose that $(K,v)$ is maximally complete of characteristic $p>0$,
and that $\F_q$ is contained in the residue field.
Then $\F_q$ can be embedded in $K$.
\end{prop}

\begin{proof}
Apply a suitable version of Hensel's lemma
(combine Theorems~6 and~7 of Chapter~2 of~\cite{schilling})
to the factorization of $x^q-x$ over $\F_q$.
\end{proof}

\begin{prop}
\label{finiteextensions}
Suppose that $(K,v)$ is maximally complete of characteristic $p>0$,
with divisible value group $G$ and with residue field $\F_q$.
If $L \subset \Kbar$ is a finite extension of $K$ of degree $n$, then
$L$ is the compositum of $\F_{q^n}$ and $K$ in $\Kbar$.
\end{prop}

\begin{proof}
Extend $v$ to $L$.
Theorem~11 in Chapter~2 of~\cite{schilling} shows that
$L$ is maximally complete,
and that $ef=n$ holds for $L/K$.
Since $G$ is divisible, there are no ordered abelian groups $G'$
with $1<(G':G)<\infty$.
Hence $e=1$, $f=n$, and the residue field of $L$ is $\F_{q^n}$.
Proposition~\ref{embedFq} implies that
the subfield $\F_{q^n}$ of $\Kbar$ is contained in $L$.
But the compositum of the linearly disjoint 
fields $\F_{q^n}$ and $K$ in $\Kbar$ is already
$n$-dimensional over $K$, so the compositum must equal $L$.
\end{proof}

\begin{rem}
	Lenstra notes that if one is interested in proving
	Theorem~\ref{degreetheorem}
	only for polynomials over $K_0=\F_q((T))$, then one can circumvent
	the theory of maximally complete fields by choosing
	$\sigma \in \Gal(\Kbar_0/K_0)$ that acts as $x \mapsto x^q$
	on $\Fbar_q$, and by taking $K$ to be the fixed field of $\sigma$.
	This $K$ contains $K_0$, still has residue field $\F_q$, and satisfies
	the conclusion of Proposition~\ref{finiteextensions}.
\end{rem}

\section{Proof of Theorem~\ref{degreetheorem}}

In proving Theorem~\ref{degreetheorem},
we may first apply Propositions~\ref{makedivisible} and~\ref{krullresult}
to assume that the value group $G$ is divisible and that
$(K,v)$ is maximally complete (still with residue field $\F_q$).
Let $\F=\F_{q^{d!}} \subset \Kbar$.
Proposition~\ref{finiteextensions}
shows that all elements of $\Kbar$ of degree at most $d$ over $K$
lie inside the compositum $L:=\F \cdot K$ of fields in $\Kbar$.
Extend $v$ to $L$.

For each $g \in G$, choose $\beta_g \in K$ with $v(\beta_g)=g$.
Now suppose $\Dbar:=\Dbar(x_0,g)$ is a closed ball in $L$.
Let $I$ be the subgroup of $\Gal(L/K) \isom \Gal(\F/\F_q)$
that maps $\Dbar$ into $\Dbar$.
Division by $\beta_g$ induces an isomorphism of $I$-modules
$\Dbar(0,g)/D(0,g) \isom \F$,
so the cohomology group $H^1(I,\Dbar(0,g)/D(0,g))$ is trivial.
The long exact sequence associated with the
exact sequence
	$$0 \rightarrow \frac{\Dbar(0,g)}{D(0,g)}
		\rightarrow \frac{L}{D(0,g)}
		\rightarrow \frac{L}{\Dbar(0,g)}
		\rightarrow 0$$
of $I$-modules shows that $\Dbar$
contains an open disk $D(x_1,g)$ mapped to itself by $I$.
We then have a bijection of $I$-sets
$\phi_{\Dbar}: \Dbar/D(0,g) \rightarrow \F$
that maps the coset $y+D(0,g)$ to the residue class of $(y-x_1)/\beta_g$.
We assume that the elements $\beta_g$ and the maps $\phi_{\Dbar}$
are fixed once and for all.

Now let $g_u$, $r$, and $S$ be as in Lemma~\ref{distances},
and let $\TT$ be the tree associated to $S$ as in Section~\ref{proof1},
so that $\ell(\TT) \le k-u$.
We now describe a labelling of the vertices of $\TT$ by elements of $\F$.
Recall that if $S_0 \in \TT$ is not a leaf,
and if $\Dbar=\Dbar(s,g)$ is the smallest disk containing $S_0$,
then the children of $S_0$ are nonempty sets of
the form $S \cap D(x_0,g)$ for some $x_0 \in \Dbar$.
Label each child by $\phi_{\Dbar}(D(x_0,g))$.
Note that the children of $S_0$ are labelled with {\em distinct}
elements of $\F$.
Finally, label the root of $\TT$
with the residue of $r/\beta_{g_u}$ in $\F^\ast$.

Let $R$ be the set of all roots of $f$ in $L$,
and let $\F[X]_{<k}$ denote the set of polynomials
of the form $a_0 + a_1 X + \cdots + a_{k-1} X^{k-1}$ with $a_i \in \F$.
We now define a map $\Phi: R \rightarrow \F[X]_{<k}$.
First, if $0 \in R$, define $\Phi(0)=0 \in \F[X]_{<k}$.
If $z \in R$ is nonzero, then $v(z)=g_u$ for some $u$.
Let $T_0 > T_1 > \cdots > T_n$ be the maximal chain ending at $T_n=\{z\}$
in the tree $\TT$ associated to $S:=R \cap D(z,g_u)$.
Define $\Phi(z)=X^{u-1} \sum_{i=0}^n \text{label}(T_i) X^i.$
Since $n \le \ell(\TT) \le k-u$, we have $\Phi(z) \in \F[X]_{<k}$.

\begin{lemma}
\label{Philemma}
$\nichts$
\begin{enumerate}
\item The map $\Phi: R \rightarrow \F[X]_{<k}$ is injective.
\item If $z \in R$ is of degree $j$ over $K$, then $\Phi(z) \in \F_{q^j}[X]$.
\end{enumerate}
\end{lemma}

\begin{proof}
To prove injectivity, we describe how to reconstruct $z$ from $\Phi(z)$.
If $\Phi(z)=0$, then $z$ must be $0$.
Otherwise its lowest degree monomial involves $X^{u-1}$ where $v(z)=g_u$.
Hence, assuming from now on that $z \not=0$,
we can reconstruct $v(z)$ from $\Phi(z)$.
Next, the {\em coefficient} of $X^{u-1}$ determines which
(nontrivial) coset of $D(0,g_u)$ in $\Dbar(0,g_u)$ $z$ belongs to.
The other coefficients uniquely determine a path ending at
the leaf $\{z\}$ in the tree associated to this coset.
Thus $\Phi(z)$ determines $z$.

For the second part, it suffices to show that
if $H$ is the subgroup of $\Gal(L/K)$ fixing $z \in R$,
then $H$ (or equivalently the isomorphic subgroup of $\Gal(\F/\F_q)$)
fixes the coefficients of $\Phi(z)$ also.
We may assume $z \not=0$.
Let $g_u=v(z)$,
and let $T_0 > T_1> \cdots > T_n=\{z\}$ be the maximal chain
in the tree $\TT$ associated to the coset $z+D(0,g_u)$ in which $z$ lies.
Since $H$ preserves the coset $z+D(0,g_u)$,
$H$ fixes the label of $T_0$.
Now suppose $1 \le i \le n$.
The smallest disk containing $T_{i-1}$ is of the form $\Dbar:=\Dbar(z,g)$
for some $g > g_u$, so $H$ is contained
in the subgroup $I \subseteq \Gal(L/K)$ preserving this disk.
The label of the child $T_i$ is $\phi_{\Dbar}(z+D(0,g))$,
and $\phi_{\Dbar}$ respects the action of $H \subseteq I$,
so $H$ fixes this label.
This holds for all $i$, so $H$ fixes all coefficients of $\Phi(z)$.
\end{proof}

Lemma~\ref{Philemma} shows that the number of zeros of $f$ in $\Kbar$
of degree at most $d$ is less than or equal to the number of polynomials
in $\F[X]_{<k}$ that are defined over $\F_{q^j}$ for some $j \le k$.
The number of such polynomials defined over $\F_{q^j}$
but no subfield is  $\sum_{i|j} q^{ik} \mu(j/i)$, by M\"obius inversion.
Theorem~2 follows upon summing over $j$.

\section*{Acknowledgements}
I thank Hendrik Lenstra for bringing the problems to my attention,
and for suggesting that $q^k$ might be the correct bound
for Theorem~\ref{maintheorem}.

\end{document}